\documentclass[12pt,reqno]{amsart}

\newtheorem{thm}{Theorem}[section]
\newtheorem{lma}[thm]{Lemma}
\newtheorem{prp}[thm]{Proposition}

\newtheorem{rem}[thm]{Remark}

\newcommand{\R}{{\mathbb R}}

\newcommand{\N}{{\mathbb N}}

\newcommand{\abs}[1]{\lvert #1 \rvert}

\newcommand{\eps}{{\varepsilon}}
\newcommand{\pf}{{\noindent \bf Proof. }}

\newcommand{\supp}{{\rm supp\ }}
\newcommand{\lb}{\langle}
\newcommand{\rb}{\rangle}

\newcommand{\dd}{\partial}

\newcommand{\ph}{\varphi}
\newcommand{\eproof}{{\mbox{\ }~\hfill
\mbox{\large $\Box$} \par \vskip 10pt}}
\newcommand{\etilde}{\,\tilde{\rule{0pt}{6pt}}\,}

\renewcommand{\d}{\partial}

\numberwithin{equation}{section}
\begin{document}

\title[Complex spherical waves in unbounded domains]{Complex spherical waves and inverse problems in unbounded domains}

\author[Salo]{Mikko Salo}
\address{Department of Mathematics and Statistics/RNI, University of Helsinki, P.O. Box 68, 00014 University of Helsinki, Finland}
\email{mikko.salo@helsinki.fi}
\author[Wang]{Jenn-Nan Wang}
\address{Department of Mathematics and NCTS Taipei, National Taiwan University,
Taipei 106, Taiwan}
\email{jnwang@math.ntu.edu.tw}

%
\begin{abstract}
This work is motivated by the inverse conductivity problem of identifying an embedded object in an infinite slab. The novelty of our approach is that we use \emph{complex spherical waves} rather than classical Calder\'on type functions. For Calder\'on type functions, they are growing exponentially on one side of   a \emph{hyperplane} and decaying exponentially on the other side. Without extra modifications, they are inadequate for treating inverse problems in  unbounded domains such as the infinite slab. The obvious reason for this is that Calder\'on type functions are not integrable on hyperplanes. So they can not be used as measurements on infinite boundaries.

For complex spherical waves used here, they blow up faster than any given positive polynomial order on the inner side of the unit sphere and decay to zero faster than any given negative polynomial order on the outer side of the unit sphere. We shall construct these special solutions for the conductivity equation in the unbounded domain by a Carleman estimate. Using complex spherical waves, we can treat the inverse problem of determining the object in the infinite slab like the problem in the bounded domain. Most importantly, we can easily localize the boundary measurement, which is of great value in practice. On the other hand, since the probing fronts are spheres, it is possible to detect some concave parts of the object.
\end{abstract}

\maketitle

\section{Introduction}

Special solutions for elliptic equations or systems have played
an important role in inverse problems since the pioneering work of
Calder\'on \cite{ca}. In 1987, Sylvester and Uhlmann \cite{syuh} introduced complex geometrical optics solutions to solve the inverse boundary
value problem for the conductivity equation. Other
type of special solutions, called oscillating-decaying solutions,
were constructed for general elliptic systems in \cite{nuwosc} and
\cite{semiosc}. These oscillating-decaying solutions have been
used in solving inverse problems, in particular in detecting
inclusions and cavities \cite{nuwosc}.

In developing the theory for inverse boundary value problems with partial
measurements, approximate complex geometrical optics solutions concentrated
near hyperplanes and near hemispheres for the Schr-\\
\"odinger
equation were given in \cite{geuh} and \cite{isuh}, respectively.
In \cite{isuh}, the construction was based on hyperbolic geometry
and was applied in \cite{iinsu} to construct complex geometrical optics
solutions for the Schr\"odinger equation where the real part of
the phase function is a radial function, i.e., its level surfaces are
spheres. These solutions are called \emph{complex spherical waves} in \cite{iinsu}. An important feature of these solutions is that they decay exponentially on one side of the sphere and grow exponentially on the other side. The hyperbolic geometry approach does not seem to work for the Laplacian with first order
perturbations, such as the Schr\"odinger equation with magnetic
potential and the isotropic elasticity. Recently, using Carleman estimates,  complex spherical waves were constructed for the Schr\"odinger
equation in \cite{ksu}, for the Schr\"odinger equation with
magnetic potential in \cite{fksu}, and for the isotropic elasticity in \cite{uhwa}. The Carleman estimate is a more flexible tool in
treating lower order perturbations.

With these complex spherical waves at hand, one can study the inverse problem of reconstructing unknown inclusions or cavities embedded in a body
with known background medium. For the conductivity equation, this problem was documented in \cite{iinsu}. There are several results, both theoretical and numerical, concerning the object identification problem by boundary measurements for the conductivity equation. We will not try to give a full account of these developments here. For detailed references, we refer to \cite{iinsu}. Here we want to point out that there are two striking advantages in using complex spherical waves for object identification problems. On one hand, one can avoid using a unique continuation procedure, more precisely, the Runge approximation, in the reconstruction. On the other hand, by the decaying property of complex spherical waves across the spherical front, the measurements can be easily localized.

The construction of complex spherical waves and the inverse problems studied in the literature cited above were restricted to bounded domains. However, there are many interesting inverse problems whose background domains are unbounded, for instance, identifying an object embedded in an infinite slab. In fact, our work here is motivated by this inverse problem. Ikehata in \cite{ikslab} studied the inverse conductivity problem in an infinite slab where the location of an inclusion is reconstructed by infinitely many boundary measurements. In \cite{ikslab}, he used Calder\'on type harmonic functions, i.e., $e^{x\cdot(\omega+i\omega^{\perp})}$ with $\omega\in{\mathbb S}^{n-1}$. These functions are not integrable on hyperplanes. Therefore, they can not be the Dirichlet data of solutions with finite energy. To remedy this, he introduced Yarmukhamedov's Green function to construct a sequence of harmonic functions with finite energy that approximate the Calder\'on type function on a bounded part of the slab and are arbitrarily small on an unbounded part of the slab.

In this paper we use complex spherical waves rather than Calder\'on type functions for the inverse problem in the slab. The most obvious advantage is that we do not need Yarmukhamedov's Green function to "localize"  complex spherical waves since we can make these solutions decay faster than any given polynomial order on infinite hyperplanes. Therefore, we can treat the problem in the infinite slab like that in the bounded domain. Furthermore, we can handle the inhomogeneous background medium without requiring the Runge approximation property. Also, since our probing fronts are spheres, we are able to determine some concave parts of the embedded object.

This paper is organized as follows. In Section~2, we discuss a Carleman estimate needed in our proof and its consequence. In Section~3, we construct complex spherical waves in the unbounded domain by means of the Carleman estimate. The investigation of the inverse problem in the infinite slab will be discussed in Section~4.

\section{Carleman estimate and related consequence}\label{sec2}

As in \cite{fksu}, \cite{ksu}, and \cite{uhwa}, we construct complex spherical waves in the unbounded domain via a Carleman estimate. Let $p\in\R^n$, $n\ge 2$, be a fixed point. Let $B_{\eps}(p)$ be the ball of radius $\eps>0$ centered at $p$. Without loss of generality, we take $p=0$. We denote $U=\overline{B_{\eps}(0)}^c$. From \cite{fksu} and \cite{ksu}, we see that $\ph(x)=\log|x|$ is a limiting Carleman weight for the semiclassical Laplacian $-h^2\Delta$ in $U$. Namely, if $a(x,\xi)=|\xi|^2-|\ph'_x|^2$ and $b(x,\xi)=2\ph'_x\cdot\xi$, then
$$\{a,b\}(x,\xi)=0\quad\text{when}\quad a(x,\xi)=b(x,\xi)=0.$$ Here $\{a,b\}=a_{\xi}'\cdot b'_x-a_x'\cdot b'_{\xi}$ is the Poisson bracket of $a, b$. With such $\ph$, we have $$e^{-\ph/h}=|x|^{-1/h}.$$ So $e^{-\ph/h}$ becomes a polynomial weight. Carleman estimates with polynomial weights are well known in proving the strong unique continuation property for Schr\"odinger operators with singular potentials, see, for example, \cite{kj}. Since we will work in $L^2$-based spaces, it is enough to use a simpler $L^2$ Carleman estimate given for instance in \cite{reg1}.
\begin{lma}\cite[Lemma~2.1]{reg1}
There exists a $C>0$ such that for any $v\in C_c^{\infty}(\R^n\setminus\{0\})$ and for any $\frac 1h\in\{k+\frac 12:k\in\N\}$ with $h\ll 1$, we have that
\begin{equation}\label{care}
{h^2}\int |x|^{-2/h}|v|^2|x|^{-n}dx\le C\int |x|^{-2/h+4}|h^2\Delta v|^2|x|^{-n}dx.
\end{equation}
Here the constant $C$ is independent of $v$ and $h$.
\end{lma}

Now we replace $1/h$ by $1/h-n-\delta$ such that $$h\in S:=\{(k+n+\delta+\frac 12)^{-1}:k\in\N\},$$ where $\delta>0$. Then \eqref{care} implies that for $v\in C_c^{\infty}(U)$ and sufficiently small $h\in S$ we have
\begin{equation}\label{cares}
{h^2}\int |x|^{-2/h}|v|^2|x|^{n+2\delta}dx\le C\int |x|^{-2/h}|h^2\Delta v|^2|x|^{n+4+2\delta}dx
\end{equation}
with possibly different constant $C$. Let $q\in L_c^{\infty}(U)$ and $P_h^{\ast}=|x|^{-1/h}h^2(\Delta-\bar q)|x|^{1/h}$. Then we get from \eqref{cares} that
\begin{equation}\label{care1}
h\|v\|_{L^2_{\frac n2+\delta}(U)}\le C\|P_h^{\ast}v\|_{L^2_{\frac n2+2+\delta}(U)}
\end{equation}
for all $v\in C_c^{\infty}(U)$ and small $h\in S$, where $L^2_s(U)$ is the weighted $L^2$ space with norm $\|f\|_{L^2_s(U)}=\|\lb x\rb^sf\|_{L^2(\Omega)}$, where $\lb x\rb=\sqrt{1+|x|^2}$. Combining \eqref{care1} and the Hahn-Banach theorem, we prove the following existence theorem.
\begin{prp}\label{exist}
For any $f\in L^2_{-\frac n2-\delta}(U)$, there exists  $w\in L^2_{-\frac n2-2-\delta}(U)$ such that
$$P_hw:=|x|^{1/h}h^2(\Delta-q)(|x|^{-1/h}w)=f\quad\text{in}\quad U$$ and
\begin{equation*}
h\|w\|_{L^2_{-\frac n2-2-\delta}(U)}\le C\|f\|_{L^2_{-\frac n2-\delta}(U)}.
\end{equation*}
\end{prp}
\pf Let us denote ${\mathcal D}=C_c^{\infty}(U)$. For $f\in L^2_{-\frac n2-\delta}(U)$, we define a linear functional $\ell$ on the linear subspace $L=P_h^{\ast}{\mathcal D}$ of $L^2_{\frac n2+2+\delta}(U)$ as follows:
$$\ell(P_h^{\ast}v)=\langle v,f\rangle_{L^2(U)}.$$ Using \eqref{care1}, we obtain
\begin{equation*}
|\ell(P_h^{\ast}v)|\le \|v\|_{L^2_{\frac n2+\delta}(U)}\|f\|_{L^2_{-\frac n2-\delta}(U)}\le \frac Ch\|P_h^{\ast}v\|_{L^2_{\frac n2+2+\delta}(U)}\|f\|_{L^2_{-\frac n2-\delta}(U)}.
\end{equation*}
So by the Hahn-Banach theorem, $\ell$ can be extended to a linear functional on $L^2_{\frac n2+2+\delta}(U)$ with the same norm. Therefore, there exists a $w\in {L^2_{-\frac n2-2-\delta}(U)}$ such that for $v\in C_c^{\infty}(\Omega)$
$$
\langle P_h^{\ast}v,w\rangle_{L^2(U)}=\langle v,f\rangle_{L^2(U)}
$$
and
$$h\|w\|_{L^2_{-\frac n2-2-\delta}(U)}\le C\|f\|_{L^2_{-\frac n2-\delta}(U)}.$$ The proof of proposition is now complete.\eproof

\section{Complex spherical waves in unbounded domains}\label{sec3}

In this section we will construct complex spherical waves for the
conductivity equation in the unbounded domain $\Omega = \{x \in
\R^n : x_n < d \}$ where $d < 0$. It is clear that
$\overline{\Omega}\subset U (=\overline{B_{\eps}(0)}^c)$ for some $\eps>0$. Let
$0<\gamma(x)\in C^{2}(\overline{\Omega})$ and $\gamma=1$ in
$B_R(0)^c\cap\Omega$ for some $R>0$. We look for a function
$v=|x|^{-1/h}\tilde v$ with $h\in S$ satisfying
\begin{equation}\label{conduc}
Lv:=\nabla\cdot(\gamma\nabla v)=0\quad\text{in}\quad \Omega.
\end{equation}
Using the Liouville transform $w=\sqrt{\gamma}v$, we know that if $w$ solves $$Pw:=-\Delta w+qw=0$$ with $q=\Delta\sqrt{\gamma}/\sqrt{\gamma}$ then $v$ solves \eqref{conduc}. Here we note that $q\in L_c^{\infty}(\Omega)$. We want to point out that by checking the proofs in \cite{iinsu}, complex spherical waves constructed there can be extended to the unbounded domain $\Omega$. In this work we take a different route using the Carleman estimate which can be adapted to other equations more easily.

As in \cite{fksu} and \cite{ksu}, we want to construct solutions to $Pw=0$ in $\Omega$ which have the form
\begin{equation*}
w = e^{-\rho(x)/h}(a+r),
\end{equation*}
where $\rho = \varphi + i\psi$ and $\varphi$, $\psi$ are given by
\begin{eqnarray*}
 & \varphi(x) = \log\,\abs{x}, & \\
 & \psi(x) = \text{dist}_{{\mathbb S}^{n-1}}(\frac{x}{\abs{x}},e_1). &
\end{eqnarray*}

Our first goal is to see what $a$ will look like. To this end, we
change coordinates and write $x = (x_1,x')$ where $x' = r\theta$
with $r > 0$ and $\theta \in {\mathbb S}^{n-2}$. We take $z = x_1 + ir$ to
be a complex variable, and let $\Psi: x \mapsto (z,\theta)$ be the
corresponding change of coordinates in $\Omega$. If $f$ is a
function in $\Omega$ we write $\tilde{f} = f \circ \Psi^{-1}$.
Thus we can see that $$\tilde{\rho} = \log\,z,$$ $$(\nabla
\rho)\etilde = \frac{1}{z}(e_1+i e_r)\quad\text{with}\quad e_r =
(0,\theta),$$ and $$(\Delta \rho)\etilde =
-\frac{2(n-2)}{z(z-\bar{z})}.$$ Also, $\nabla \rho \cdot \nabla$
becomes $(2/z)\partial_{\bar{z}}$ in the new coordinates.

Writing $h^2 P = (-ih\nabla)^2 + h^2 q$, we have
\begin{eqnarray}\label{comp}
e^{\rho/h} h^2 P e^{-\rho/h}&=& (-ih\nabla+i\nabla \rho)^2 + h^2 q\notag\\
&=& -(\nabla \rho)^2 + h(2\nabla\rho \cdot \nabla + \Delta \rho) + h^2 P\notag\\
&=& h(2\nabla\rho \cdot \nabla + \Delta \rho) + h^2 P.
\end{eqnarray}
Here, $\varphi$ and $\psi$ were chosen so that $(\nabla \rho)^2 = 0$. To get $P(e^{-\rho/h}(a+r)) = 0$, we choose $a$ and $r$ to satisfy in $\Omega$
\begin{equation}\label{transp}
(\nabla \rho \cdot \nabla + \frac{1}{2} \Delta \rho)a = 0
\end{equation}
and
\begin{equation}\label{r}
e^{\rho/h} h^2 P e^{-\rho/h} r = -h^2 Pa.
\end{equation}
The first equation \eqref{transp} is a transport equation for $a$. Writing $a = \tilde{a} \circ \Psi$, this reduces to
\begin{equation*}
(\partial_{\bar{z}} - \frac{n-2}{2(z-\bar{z})}) \tilde{a} = 0 \quad \text{in } \Psi(\Omega).
\end{equation*}
It is easy to see that the general solution is given by
\begin{equation*}
\tilde{a} = (z-\bar{z})^{\frac{2-n}{2}} \tilde{g}
\end{equation*}
for any $\tilde{g}(z,\theta)$ satisfying $\partial_{\bar{z}} \tilde{g} = 0$.

Taking $\tilde{g} = 1$ and going back to the $x$ coordinates, we obtain
\begin{equation}\label{a}
a(x) = (2i\abs{x'})^{\frac{2-n}{2}}, \quad x \in \Omega.
\end{equation}
In \eqref{r}, the right hand side becomes
\begin{equation*}
-h^2 Pa = (2i)^{\frac{2-n}{2}} h^2 \Big( \frac{1}{4} (n-2)(4-n) \abs{x'}^{\frac{-2-n}{2}} - q \abs{x'}^{\frac{2-n}{2}} \Big).
\end{equation*}
This decays in $x'$ but is constant in $x_1$, so $Pa \in L^2_{s}(\Omega)$ for $s < -n/2$. We now rewrite \eqref{r} as
$$
-|x|^{1/h}h^2P(|x|^{-1/h}e^{-i\psi/h}r)=e^{-i\psi/h}h^2Pa.
$$
Since $e^{-i\psi/h}h^2Pa\in L^2_{s}(\Omega)$ with $s=-n/2-\delta$ for some $\delta>0$, we deduce from Proposition~\ref{exist} that there exists $r$ solving \eqref{r} and satisfying
\begin{equation}\label{est1}
\|r\|_{L^2_{-\frac n2-2-\delta}(\Omega)}\le Ch.
\end{equation}
Therefore, we have constructed the special solution $w$ of
$-\Delta w+qw=0$ in $\Omega$ having the form
$$w=e^{-(\ph+i\psi)/h}(a+r)=|x|^{-1/h}e^{-i\psi/h}(a+r)$$
with $a$ given by \eqref{a} and $r$ satisfying \eqref{est1}. In
particular, we can see that $w\in L^2(\Omega)$ for all
sufficiently small $h$. Furthermore, since $\Delta w=qw\in
L^2(\Omega)$, we can use a cut-off technique and the elliptic estimate to show that
$$w\in H^2(\tilde\Omega),$$ where $\tilde{\Omega}=\{x\in\R^n :
x_n<\tilde d\}$ for any $\tilde d<d$. In other words, the trace of
$w$ on any hyperplane $H_b:=\{x_n=b\}$ with $b<\tilde d$ is
well-defined and
$$\d_{x_n}w|_{H_b}\in H^{1/2}.$$ This property is an important difference
between complex spherical waves and Calder\'on type solutions.

In order to apply complex spherical waves to the inverse problem, we need to estimate the $H^1$ norm of $r$ on any fixed open subset of $\Omega$. Let $\Omega_1$ and $\Omega_2$ be two open bounded subsets of $\Omega$ such that $\overline{\Omega_1}\subset\Omega_2$ and $\overline{\Omega_2}\subset\Omega$. From \eqref{comp} we see that
$$e^{\rho/h} h^2 P e^{-\rho/h}=-h^2\Delta+h(2\nabla\rho \cdot \nabla + \Delta \rho) + h^2 q.$$ Therefore, $e^{\rho/h} h^2 P e^{-\rho/h}$ is elliptic as a semiclassical operator. Furthermore, we have
$$\|e^{\rho/h} h^2 P e^{-\rho/h}r\|_{L^2(\Omega_2)}=\|h^2Pa\|_{L^2(\Omega_2)}\le Ch^2$$ and $$\|r\|_{L^2(\Omega_2)}\le Ch.$$ Hence, by the semiclassical version of elliptic estimate \cite[Lemma~2.6]{ktz}, we obtain that
$$\sum_{|\alpha|\le 2}\|h^{\alpha}\d^{\alpha}r\|_{L^2(\Omega_1)}\le Ch,$$ in particular,
\begin{equation}\label{est22}
\|r\|_{L^2(\Omega_1)}+h\|\nabla r\|_{L^2(\Omega_1)}\le Ch.
\end{equation}
Finally, let $p$ be any point such that $\text{dist}(p,\Omega)>0$, then complex spherical waves for \eqref{conduc} in $\Omega$ are given by
$$v=v(x,h)=\gamma^{-1/2}|x-p|^{-1/h}e^{-i\psi/h}(a+r)$$ for $h\in S$ small, where $\psi=\text{dist}_{{\mathbb S}^{n-1}}(\frac{x-p}{\abs{x-p}},e_1)$, $a(x)=(2i\abs{x'-p'})^{\frac{2-n}{2}}$ and $r$ satisfies \eqref{est1} and \eqref{est22}. Note that for all $x\in\Omega$ we have $x'-p'=(x_2-p_2,\cdots,x_n-p_n)\neq 0$ since $p_n>x_n$.

\section{Inverse problem in a slab}\label{sec4}

In this section we shall apply the special solutions we
constructed above to the inverse problem in an infinite slab. We
study the reconstruction of an embedded cavity here. The same
method works for the inclusion case and so does the method in
\cite{iinsu} or \cite{uhwa}. We leave this generalization to the
interested reader. Let $\Omega=\{x\in\R^n : d_1<x_n<d_2,\
d_1<d_2\}.$ Assume that $D$ is a domain with $C^2$ boundary in
$\Omega$ so that $\bar D\subset\Omega$ and $\Omega\setminus\bar D$
is connected. Let $u(x)$ be the unique solution of finite energy
to
\begin{equation}\label{meq}
\begin{cases}
L_{\gamma}u=0\quad\text{in}\quad\Omega\setminus\bar D,\\
\gamma\frac{\partial u}{\partial\nu}=0\quad\text{on}\quad\partial D,\\
u=f\in H^{1/2}\quad\text{on}\quad\partial\Omega,
\end{cases}
\end{equation}
where the conductivity parameter $\gamma(x)$ satisfies the assumptions described in the previous section. The well-posedness of \eqref{meq} can be proved by the standard Lax-Milgram theorem with a Poincar\'e-type inequality in the infinite slab: there exists $C>0$ independent of $u$ such that
\begin{equation}\label{poincare}
\int_{\Omega\setminus\bar D}|u|^2dx\le C\int_{\Omega\setminus\bar D}|\nabla u|dx
\end{equation}
for all $u\in H^1(\Omega\setminus\bar D)$ with $u=0$ on $x_n=d_1$ and $d_2$. One possible way to prove \eqref{poincare} is to divide $\Omega= Z\cup(\Omega\setminus Z)$, where $Z$ is a cylindrical domain of the form $Z = B \times (d_1,d_2)$, $B$ is a ball in $\R^{n-1}$, and $\bar D\subset Z$. Then \eqref{poincare} is a consequence of
\begin{equation}\label{poin1}
\int_{Z}|v|^2dx\le C\int_{Z}|\nabla v|dx
\end{equation}
for all $v\in H^1(Z)$ with $v=0$ on $\partial Z\cap\{x_n=d_1, d_2\}$ and
\begin{equation}\label{poin2}
\int_{\Omega\setminus\bar Z}|w|^2dx\le C\int_{\Omega\setminus\bar Z}|\nabla w|dx
\end{equation}
for all $w\in H^1(\Omega\setminus\bar Z)$ with $w=0$ on $\partial (\Omega\setminus\bar Z)\cap\{x_n=d_1, d_2\}$. Now \eqref{poin1} can be proved by a contradiction argument and \eqref{poin2} can be established by the usual integration technique. For brevity, we left the details to the reader.

The inverse problem here is to identify $D$ from the Dirichlet-to-Neumann map $\Lambda_D:H^{1/2}(\partial\Omega)\mapsto H^{-1/2}(\partial\Omega)$ $$\Lambda_D:f\mapsto \gamma\frac{\partial u}{\partial\nu}|_{\partial\Omega}.$$ In the weak formulation, the Dirichlet-to-Neumann map is defined by
$$\lb\Lambda_Df,g\rb=\int_{\Omega\setminus\bar D}\gamma\nabla u\cdot\nabla vdx,$$ where $v\in H^1(\Omega\setminus\bar D)$ with $g=v|_{\partial\Omega}$. We are interested in the reconstruction problem in this work.

Let $\hat\Omega=\{x\in\R^3: x_n<\hat d\}$ with $d_2<\hat d$ and $t>0$. Pick any point $p$ so that $\text{dist}(p,\hat\Omega)>0$. We now denote complex spherical waves of $L_{\gamma}v=0$ in $\hat\Omega$ by
$$v_t(x,h)=t^{1/h}\gamma^{-1/2}|x-p|^{-1/h}e^{-i\psi/h}(a+r)=\gamma^{-1/2}(\frac{t}{|x-p|})^{1/h}e^{-i\psi/h}(a+r).$$ Let us define the energy gap functional
$$E_t(h):=\langle(\Lambda_0-\Lambda_D)\bar v_t(x,h),v_t(x,h)\rangle=\int_{\partial\Omega}(\Lambda_0-\Lambda_D)\bar v_t(x,h)\cdot v_t(x,h)\ ds,$$ where $\Lambda_0$ is the Dirichlet-to-Neumann map for $L_{\gamma}$ with $D=\emptyset$. Note that $E_t(h)$ is well-defined for any $t>0$ even when $\partial\Omega$ is unbounded. To understand how we reconstruct $D$, we first observe that $E_t(h)$ can be estimated by
\begin{equation}\label{est}
\frac 1C\int_D|\nabla v_t(x,h)|^2dx\le E_t(h)\le C\int_D(|\nabla v_t(x,h)|^2+|v_t(x,h)|^2)dx
\end{equation}
for some constant $C>0$.

\medskip\noindent{\bf Proof of} \eqref{est}.
Let $u$ be the solution of \eqref{meq} with $f=v_t$. We observe that
\begin{equation}\label{eq1}
\int_{\Omega\setminus\bar D}\gamma\nabla u\cdot\nabla(\bar u-\bar v_t)dx=0.
\end{equation}
Using the definition of the Dirichlet-to-Neumann map and \eqref{eq1}, we can compute
\begin{eqnarray}\label{eq2}
&&\langle(\Lambda_0-\Lambda_D)\bar v_t,v_t\rangle\notag\\
&=&\int_{\dd\Omega}\gamma\frac{\dd\bar v_t}{\dd\nu}v_tds-\int_{\dd\Omega}\gamma\frac{\dd\bar u}{\dd\nu}v_tds\notag\\
&=&\int_{\Omega}\gamma\nabla\bar v_t\cdot\nabla v_tdx-\int_{\Omega\setminus\bar D}\gamma\nabla\bar u\cdot \nabla v_tdx\notag\\
&=&\int_{\Omega\setminus\bar D}\gamma\nabla\bar v_t\cdot\nabla v_tdx-\int_{\Omega\setminus\bar D}\gamma\nabla\bar u\cdot \nabla v_tdx+\int_{\Omega\setminus\bar D}\gamma\nabla u\cdot\nabla(\bar u-\bar v_t)dx\notag\\
&&+\int_D\gamma|\nabla v_t|^2dx\notag\\
&=&\int_D\gamma|\nabla v_t|^2dx+\int_{\Omega\setminus\bar D}\gamma\nabla(u-v_t)\cdot\nabla(\bar u-\bar v_t)dx.
\end{eqnarray}
The formula \eqref{eq2} immediately implies the first inequality of \eqref{est}.

To obtain the second inequality of \eqref{est}, we need to estimate the last term in \eqref{eq2}. Denote $w=u-v_t$. By the construction, we have that
$$
\begin{cases}
L_{\gamma} w=0\quad\text{in}\quad\Omega\setminus\bar D,\\
\gamma\frac{\partial w}{\partial\nu}=-\gamma\frac{\partial v_t}{\partial\nu}\quad\text{on}\quad\partial D,\\
w=0\quad\text{on}\quad\partial\Omega.
\end{cases}
$$
It follows from the elliptic regularity theorem that
\begin{equation}\label{eq3}
\|w\|_{H^1(\Omega\setminus\bar D)}\le C\|\gamma\frac{\partial v_t}{\partial\nu}\|_{H^{-1/2}(\dd D)}.
\end{equation}
On the other hand, we know that $L_{\gamma} v_t=0$ in $D$ and therefore,
\begin{equation}\label{eq5}
\|\gamma\frac{\partial v_t}{\partial\nu}\|_{H^{-1/2}(\dd D)}\le C\|v_t\|_{H^1(D)}.
\end{equation}
Putting together \eqref{eq3} and \eqref{eq5} leads to
$$
\|w\|_{H^1(\Omega\setminus\bar D)}\le C\|v_t\|_{H^1(D)}
$$
and the second estimate in \eqref{est}.\eproof

Now we will use \eqref{est} to study the behaviors of $E_t(h)$ as $h\to 0$ and $h\in S$ for different $t$'s.
\begin{thm}\label{full}
Let $\text{\rm dist}(D,p)=:d_0>0$. For $t>0$ and sufficiently small $h\in S$, we have that

\noindent{\rm(i).} If $d_0>t$, then
$E_t(h)\le C \alpha^{1/h}$ for some  $\alpha<1$;

\noindent{\rm(ii).} If $d_0<t$, then $E_t(h)\ge C \beta^{1/h}$ for
some $\beta>1$;

\noindent{\rm(iii).} If $\overline{D}\cap\overline{B_t(p)}=\{y\}$,
then ${C}^{-1}h^{n-2}\le E_t(h)\le Ch^{-1}$ for some $C>0$.
\end{thm}
\pf (i). If $\text{\rm dist}(D,p)=d_0>t$ then
$$\frac{t}{|x-p|}<1\quad\forall\ x\in\bar D.$$ Combining the second inequality of \eqref{est} and the behavior of $v_t(x,h)$, we immediately prove this statement.

\medskip\noindent (ii). We first pick a small ball
$B_{\delta}\subset\subset B_{t}(p)\cap D$. Using \eqref{est22}, the leading term of $\nabla v_t(x,h)$ is
\begin{equation}\label{leader}
-\frac 1h(\frac{t}{|x-p|})^{1/h}e^{-i\psi/h}\gamma^{-1/2}\nabla(\log|x-p|+i\psi)a.
\end{equation}
We notice that $\gamma^{-1/2}\nabla(\log|x-p|+i\psi)a\ne 0$ for all $x\in B_{\delta}$. So this statement follows from the first inequality of \eqref{est} and the fact
$$\frac{t}{|x-p|}>1\quad\forall\ x\in B_{\delta}.$$

\medskip\noindent (iii). Pick a small cone with vertex at $y$, say $\Gamma$, so that there
exists an $\epsilon>0$ satisfying
$$
\Gamma\cap\{0<|x-y|<\epsilon\}\subset D.
$$
It is not restrictive to take $y=0$ and $p=(0,\cdots,0,t)$ for $t>0$. Now we observe that if $z\in\Gamma$ and $|z-y|=|z|=s<\epsilon$ then
$$|z-p|\le s+t,$$ that is $$\frac{t}{|z-p|}\ge\frac{t}{s+t}.$$ Hence, from the first inequality of \eqref{est} and \eqref{leader}, for $h\ll 1$ we have that
\begin{eqnarray*}
E_t(h)&\ge& \frac{C}{h^2}\int_{0}^{\epsilon}(\frac{t}{s+t})^{2/h}s^{n-1}ds\\
&\ge&\frac{C}{h^2}\sum_{k=1}^n\frac{h}{2-kh}\\
&\ge&Ch^{n-2}.
\end{eqnarray*}

On the other hand, we can choose a cone $\widetilde{\Gamma}$ with
vertex at $p$ such that
$\overline{D}\subset\widetilde{\Gamma}\cap\{t\le|x-p|<t+\eta\}$ for $\eta>0$. Thus, by the second inequality of \eqref{est} and \eqref{leader},
we can estimate
\begin{eqnarray*}
E_t(h)&\le&C(t+\eta)^{n-1}\frac{1}{h^2}\int_{t}^{t+\eta}(\frac{t}{s})^{2/h}ds\\
&\le&Ch^{-1}.
\end{eqnarray*}
\eproof

The distinct behavior of $E_t(h)$ immediately allows us to detect the boundary of the cavity. However, Theorem~\ref{full} is unpractical since we need to take the measurements on the whole unbounded boundary. Fortunately, taking advantage of the decaying property of complex spherical waves, we are able to localize the measurement near the part where complex spherical waves are not decaying. To be precise, let $\phi_{\delta,t}(x)\in
C_0^{\infty}(\R^n)$ satisfy
$$
\phi_{\delta,t}(x)=\begin{cases} 1\quad\text{on}\quad B_{t+\delta/2}(p)\\
0\quad\text{on}\quad\R^n\setminus\overline{B_{t+\delta}(p)}
\end{cases}
$$
where $\delta>0$ is sufficiently small. Now we are going to use the measurement
$f_{\delta,t}(x,h)=\phi_{\delta,t}v_{t}(x,h)|_{\partial\Omega}$.
Clearly, the measurement $f_{\delta,t}$ is localized on
$B_{t+\delta}(p)\cap\partial\Omega$. In fact, if $\delta$ is sufficiently small, $f_{\delta,t}$ is localized on only one part of $\partial\Omega$, depending on whether $p$ lies above or below $\Omega$.  Let us define
$$
E_{\delta,t}(h)=\langle(\Lambda_0-\Lambda_D)\bar{f}_{\delta,t},f_{\delta,t}\rangle.
$$
\begin{thm}\label{local}
The statements of Theorem~\ref{full} are valid for
$E_{\delta,t}(h)$.
\end{thm}
\pf The main idea is to prove that the error caused by the
remaining part of the measurement
$g_{\delta,t}:=(1-\phi_{\delta,t})v_t(x,h)|_{\partial\Omega}$ is as
small as any given polynomial order. Let ${w}_{\delta,t}(x,h)$ be the unique solution of
$$L_{\gamma}w=0\quad\text{in}\quad\Omega$$
with boundary value $g_{\delta,t}$. We now want to compare
$w_{\delta,t}$ with $(1-\phi_{\delta,t})v_{t}$. To this end,
we first observe that
$$
\begin{cases}
L_{\gamma}((1-\phi_{\delta,t})v_{t}-w_{\delta,t})=L_{\gamma}((1-\phi_{\delta,t})v_{t})\quad\text{in}\quad\Omega,\\
(1-\phi_{\delta,t})v_{t}-w_{\delta,t}=0\quad\text{on}\quad\partial\Omega.
\end{cases}
$$
Notice that $$L_{\gamma}((1-\phi_{\delta,t})v_{t})=(1-\phi_{\delta,t})L_{\gamma}v_t+[L_{\gamma},(1-\phi_{\delta,t})]v_t=[L_{\gamma},(1-\phi_{\delta,t})]v_t.$$ So we see that $$\supp(L_{\gamma}((1-\phi_{\delta,t})v_{t}))\subset\overline\Omega\cap\{t+\delta/2\le |x-p|\le t+\delta\}$$ and therefore
$$
\|L_{\gamma}((1-\phi_{\delta,t})v_{t})\|_{L^2(\Omega)}\le
C{\alpha}_1^{1/h}
$$ for some $0<\alpha_1<1$. Consequently, we have that
$$
\|(1-\phi_{\delta,t})v_{t}-w_{\delta,t}\|_{H^1(\Omega)}\le
C\alpha_1^{1/h},
$$
in particular,
\begin{equation}\label{err}
\|(1-\phi_{\delta,t})v_{t}-w_{\delta,t}\|_{H^1(D)}\le
C\alpha_1^{1/h}.
\end{equation}
Using \eqref{est} for
$\langle(\Lambda_0-\Lambda_D)\bar g_{\delta,t},g_{\delta,t}\rangle$
with $v_t$ being replaced by $w_{\delta,t}$, we get from
\eqref{err} and decaying property of $v_{t}$ that
$$
\langle(\Lambda_0-\Lambda_D)\bar g_{\delta,t},g_{\delta,t}\rangle\le C{\alpha}_2^{1/h}
$$ for some $0<{\alpha}_2<1$.

Now we first consider (i) of Theorem~\ref{full} for $E_{\delta,t}(h)$. We shall use a trick given in \cite{iinsu}. In view of the definition of the energy gap functional and the proof of the first inequality of \eqref{est}, we get that
$$
0\le \langle(\Lambda_0-\Lambda_D)(\zeta \bar f_{\delta,t}\pm\zeta^{-1}\bar g_{\delta,t}),{\zeta f_{\delta,t}\pm\zeta^{-1}g_{\delta,t}}\rangle
$$
for any $\zeta>0$, which leads to
\begin{eqnarray}\label{eq200}
&&|\langle(\Lambda_0-\Lambda_D)\bar f_{\delta,t},{g_{\delta,t}}\rangle+\langle(\Lambda_0-\Lambda_D)\bar g_{\delta,t},{f_{\delta,t}}\rangle|\notag\\
&\le& \zeta^2\langle(\Lambda_0-\Lambda_D)\bar f_{\delta,t},{f_{\delta,t}}\rangle+\zeta^{-2}\langle(\Lambda_0-\Lambda_D)\bar g_{\delta,t},{g_{\delta,t}}\rangle.
\end{eqnarray}
It now follows from $v_t(x,h)|_{\partial\Omega}=f_{\delta,t}+g_{\delta,t}$ and \eqref{eq200} with $\zeta=1/\sqrt{2}$ that
\begin{eqnarray}\label{eq300}
&&\frac 12\langle(\Lambda_0-\Lambda_D)\bar f_{\delta,t},{f_{\delta,t}}\rangle\notag\\
&\le& \langle(\Lambda_0-\Lambda_D)\bar g_{\delta,t},{g_{\delta,t}}\rangle+\langle(\Lambda_0-\Lambda_D)\bar v_{t},{v_{t}}\rangle\notag\\
&\le& C{\alpha}_2^{1/h}+\langle(\Lambda_0-\Lambda_D)\bar v_{t},{v_{t}}\rangle.
\end{eqnarray}
So from (i) of Theorem~\ref{full}, the same statement holds for $E_{\delta,t}(h)$. Likewise, the second inequality of (iii) in Theorem~\ref{full} also holds.

Next we consider (ii) and the first inequality of (iii) in Theorem~\ref{full} for $E_{\delta,t}$. Choosing $\zeta=1$ in \eqref{eq200} we get that
\begin{eqnarray}\label{eq500}
&&\frac 12\lb(\Lambda_0-\Lambda_D)\bar v_{t},v_{t}\rb\notag\\
&\le& \lb(\Lambda_0-\Lambda_D)\bar g_{\delta,t},{g_{\delta,t}}\rb+\lb(\Lambda_0-\Lambda_D)\bar f_{\delta,t},{f_{\delta,t}}\rb\notag\\
&\le& C{\alpha}_2^{1/h}+\lb(\Lambda_0-\Lambda_D)\bar f_{\delta,t},{f_{\delta,t}}\rb.
\end{eqnarray}
Therefore, (ii) of Theorem~\ref{full} and \eqref{eq500} implies that the same fact holds for $E_{\delta,t}$. Likewise, the first inequality of (iii) in Theorem~\ref{full} holds  for $E_{\delta,t}$. The proof is now complete.\eproof

\begin{rem}
With the help of Theorem~\ref{local}, when parts of $\partial D$ are near the boundary $\partial\Omega$, we could detect some partial information of $\partial D$
from only a few measurements taken from a very small region on one side of $\partial\Omega$.
\end{rem}

To end the presentation, we provide an algorithm of the method.

\medskip
Step 1. Pick a point $p\notin\overline\Omega$ and near $\partial\Omega$. Construct
complex spherical waves $v_{t}(x,h)$ for $h\in S$.

\smallskip
Step 2. Draw two balls $B_t(p)$ and $B_{t+\delta}(p)$. Set the
Dirichlet data
$f_{\delta,t}=\phi_{\delta,t}v_{t}|_{\partial\Omega}$. Measure
the Neumann data $\Lambda_Df_{\delta,t}$ over the region
$B_{t+\delta}(p)\cap\partial\Omega$.

\smallskip
Step 3. Calculate
$E_{\delta,t}(h)=\lb(\Lambda_0-\Lambda_D)\bar f_{\delta,t},f_{\delta,t}\rb$.
If $E_{\delta,t}(h)$ tends to zero as $h\to 0$, then the probing front
$\{|x-p|=t\}$ does not intersect the inclusion. Increase $t$ and
compute $E_{\delta,t}(h)$ again.

\smallskip
Step 4. If $E_{\delta,t}(h)$ increases to $\infty$ as $h\to 0$,
then the front $\{|x-p|=t\}$ intersects the inclusion. Decrease
$t$ to make more accurate estimate of $\partial D$.

\smallskip
Step 5. Choose a different $p$ and repeat the procedures Step 1-4.

\section*{Acknowledgements}

Mikko Salo was partially supported by the Academy of Finland. Jenn-Nan Wang was supported in part by the National Science Council of Taiwan (NSC 94-2115-M-002-003). This work was done when both authors were visiting the University of Washington. We would like to thank Gunther Uhlmann for his
encouragements and the Department of Mathematics at the University
of Washington for its hospitality.


\end{document}